\documentclass[graybox]{svmult}

\usepackage{type1cm}        % activate if the above 3 fonts are
\usepackage{makeidx}         % allows index generation
\usepackage{graphicx}        % standard LaTeX graphics tool
\usepackage{multicol}        % used for the two-column index
\usepackage[bottom]{footmisc}% places footnotes at page bottom

\usepackage{newtxtext}       % 
\usepackage{hyperref}
\usepackage{amsmath}
\usepackage{amsfonts}
\usepackage{mathtools}
\usepackage{braket} % for \set
\DeclareMathOperator{\dive}{div}
\usepackage{subfig}

\makeindex             % used for the subject index
\begin{document}

\title*{How does the partition of unity influence SORAS preconditioner?}
\author{Marcella Bonazzoli, Xavier Claeys, Fr\'ed\'eric Nataf and Pierre-Henri Tournier}
\authorrunning{M. Bonazzoli, X. Claeys, F. Nataf and P.-H. Tournier}
\institute{Marcella Bonazzoli \at Inria, UMA, ENSTA Paris, Institut Polytechnique de Paris, Palaiseau, France, \email{marcella.bonazzoli@inria.fr}
\and Xavier Claeys \at Sorbonne Universit\'e, CNRS, Universit\'e Paris Cit\'e, LJLL, Paris, France, \email{xavier.claeys@sorbonne-universite.fr}, \and Fr\'ed\'eric Nataf \at Sorbonne Universit\'e, CNRS, Universit\'e Paris Cit\'e, Inria, LJLL, Paris, France, \email{frederic.nataf@sorbonne-universite.fr}, \and Pierre-Henri Tournier \at Sorbonne Universit\'e, CNRS, Universit\'e Paris Cit\'e, Inria, LJLL, Paris, France, \email{pierre-henri.tournier@sorbonne-universite.fr}}
\maketitle

\abstract*{We investigate the influence of the choice of the partition of unity on the convergence of the Symmetrized Optimized Restricted Additive Schwarz (SORAS) preconditioner for the reaction-convection-diffusion equation. We focus on two kinds of partitions of unity, and study the dependence on the overlap and on the number of subdomains. In particular, the second kind of partition of unity, which is non-zero in the interior of the whole overlapping region, gives more favorable convergence properties, especially when increasing the overlap width, in comparison with the first kind of partition of unity, whose gradient is zero on the subdomain interfaces and which would be the natural choice for ORAS solver instead.}

\section{Introduction}

The Symmetrized Optimized Restricted Additive Schwarz (SORAS) preconditioner, first introduced in \cite{KiSa:OBDD:2007} for the Helmholtz equation and called OBDD-H, was later studied in \cite{HaJoNa:soras:2015} for generic symmetric positive definite problems and viewed as a symmetric variant of ORAS preconditioner. 
Its convergence was rigorously analyzed in \cite{GrSpZo:impedance} for the Helmholtz equation, and in \cite{BCNT:2021:SORAS} we generalized this theory to generic non self-adjoint or indefinite problems. Moreover, as an illustration of our theory, we proved new estimates for the specific case of the heterogeneous reaction-convection-diffusion equation. In the numerical experiments in \cite{BCNT:2021:SORAS}, we noticed that the number of iterations for convergence of preconditioned GMRES appears not to vary significantly when increasing the overlap width. In the present paper, we show that actually this is due to the particular choice of the partition of unity for the preconditioner. The influence of five different kinds of partition of unity on SORAS solver and preconditioner for the Laplace equation has been briefly studied in the conclusion of \cite{Gander:2020:PUdd25}, where the method is named ORASH. Here, for the reaction-convection-diffusion equation, we focus on two kinds of  partitions of unity, and study the dependence on the overlap and on the number of subdomains.

\section{SORAS preconditioner and two kinds of partition of unity}

Let $A$ denote the $n \times n$ matrix, not necessarily positive definite nor self-adjoint, arising from the discretization of the problem to be solved, posed in an open domain $\Omega \subset \mathbb{R}^d$. 
Given a set of overlapping open subdomains $\Omega_j, j=1,\dots,N$, such that $\Omega = \cup_{j=1}^N \Omega_j$ and each $\overline{\Omega_j}$ is a union of elements of the mesh $\mathcal{T}^h$ of $\Omega$, we consider the set $\mathcal{N}$ of the unknowns on the whole domain, so $\#\mathcal{N} = n$, and its decomposition $\mathcal{N} = \bigcup_{j=1}^N \mathcal{N}_j$ into the non-disjoint subsets corresponding to the different overlapping subdomains {$\overline{\Omega}_j \cap \Omega$}, with $\#\mathcal{N}_j = n_j$. Denote by $\delta$ the width of the overlap between subdomains. The following  matrices are then the classical ingredients to define overlapping Schwarz domain decomposition preconditioners (see e.g.~\cite[\S1.3]{DoJoNa:bookDDM}):
\begin{itemize}
\item
{restriction} matrices $R_j$ from $\Omega$ to {$\overline{\Omega}_j \cap \Omega$}, which are $n_j \times n$ Boolean matrices whose $(i,i')$ entry equals $1$ if the $i$-th unknown in $\mathcal{N}_j$ is the $i'$-th one in $\mathcal{N}$ and vanishes otherwise; 
\item
{extension} by zero matrices $R^T_j$ from {$\overline{\Omega}_j \cap \Omega$} to $\Omega$;
\item
{partition of unity} matrices $D_j$, which are $n_j \times n_j$ diagonal matrices with real non-negative entries such that $\sum_{j=1}^N R_j^T D_j R_j = I$ and which can be seen as matrices that properly weight the unknowns belonging to the overlap between subdomains; 
\item
{local matrices} $B_j$, of size $n_j \times n_j$, which arise from the discretization of subproblems posed in {$\overline{\Omega}_j \cap \Omega$}, with for instance Robin-type or more general absorbing transmission conditions on the interfaces $\partial\Omega_j \setminus \partial\Omega$. 
\end{itemize}
Then the one-level Symmetrized Optimized Restricted Additive Schwarz (SORAS) preconditioner is defined as 
\begin{equation}
\label{eq:SORAS}
M^{-1} \coloneqq \sum_{j=1}^N R_j^T D_j B_j^{-1} D_j R_j. 
\end{equation}
Note that $M^{-1}$ is not self-adjoint when $B_j$ is not self-adjoint, even if we maintain the SORAS name, where S stands for `Symmetrized'. 
In fact, this denomination was introduced in \cite{HaJoNa:soras:2015} for symmetric positive definite problems, since in that case SORAS preconditioner is a symmetric variant of ORAS preconditioner $\sum_{j=1}^N R_j^T D_j B_j^{-1} R_j$. Thus, the adjective `Symmetrized' stands for the presence of the rightmost partition of unity $D_j$. We recall that `Restricted' indicates the presence of the leftmost partition of unity $D_j$ and that `Optimized' refers to the choice of transmission conditions other than standard Dirichlet conditions in the local matrices $B_j$.

Here we focus on the influence exerted by the choice of partition of unity matrices $D_j$ on the convergence of GMRES preconditioned by \eqref{eq:SORAS}. Indeed, several definitions of the diagonal matrices $D_j$ are possible to ensure property $\sum_{j=1}^N R_j^T D_j R_j = I$.  In general, the diagonals of the $D_j$ can be constructed by the interpolation of continuous partition of unity functions $\chi_j \colon \Omega \to [0,1]$, $j=1,\dots,N$: $\sum_{j=1}^N \chi_j = 1$ in $\overline{\Omega}$, and $\mathrm{supp}( \chi_j )\subset \Omega_j$, so in particular $\chi_j$ is zero on the subdomain interfaces $\partial \Omega_j \setminus \partial \Omega$. 

In addition, in the case of ORAS fixed-point iterative solver, also the first derivatives of $\chi_j$ are required to be equal to zero on $\partial \Omega_j \setminus \partial \Omega$, because this property ensures that the continuous version of ORAS solver is equivalent to Lions' algorithm, see e.g.~\cite[\S2.3.2]{DoJoNa:bookDDM} for a particular model problem. An instructive calculation for a simple one- (and two-) dimensional problem, which shows an analogous equivalence property for RAS solver, is given in \cite{EfGa:2003:RAS}; a more general equivalence result for ORAS solver is proved in \cite[Theorem 3.4]{StGaTh:2007:OMA}. 
\begin{figure}
\centering
\subfloat[][{PU1, $\delta=2h$}]
    {\includegraphics[width=.48\textwidth]{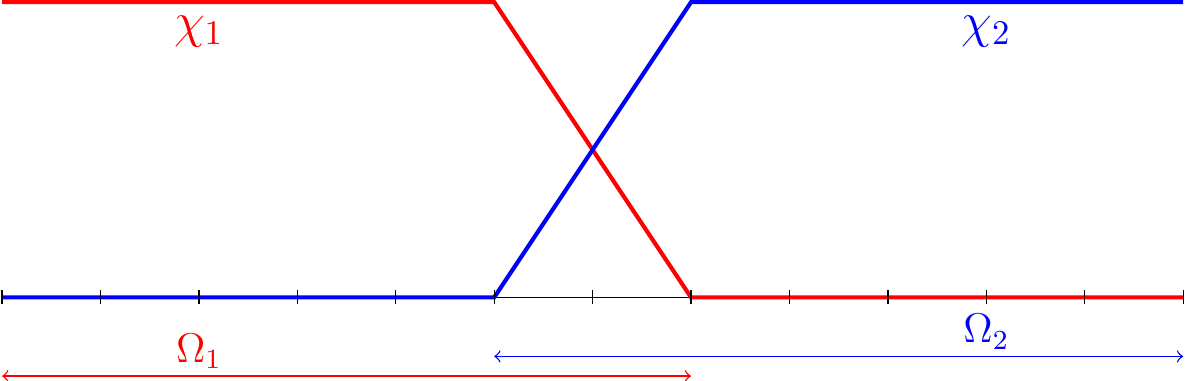}} \quad
\subfloat[][{PU2, $\delta=2h$}]
    {\includegraphics[width=.48\textwidth]{default_and_raspart_PU_ovr1.pdf}} \\
\subfloat[][{PU1, $\delta=4h$}]
    {\includegraphics[width=.48\textwidth]{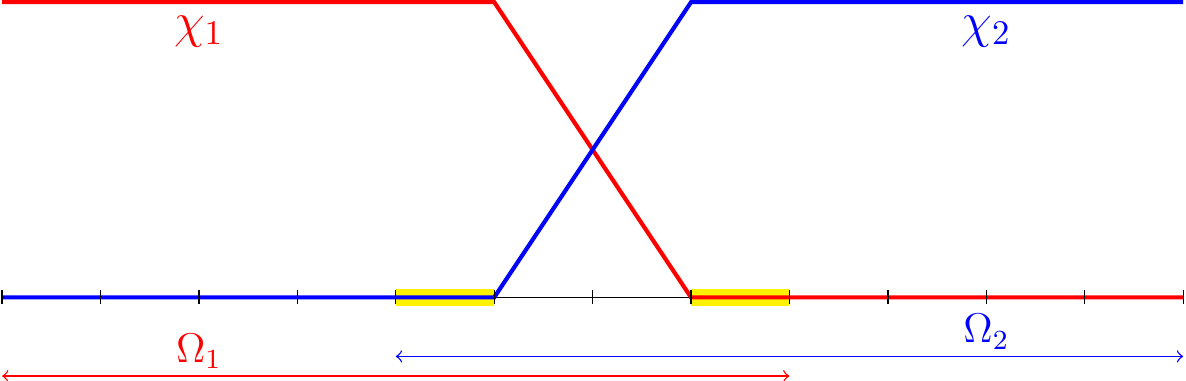}} \quad
\subfloat[][{PU2, $\delta=4h$}]
    {\includegraphics[width=.48\textwidth]{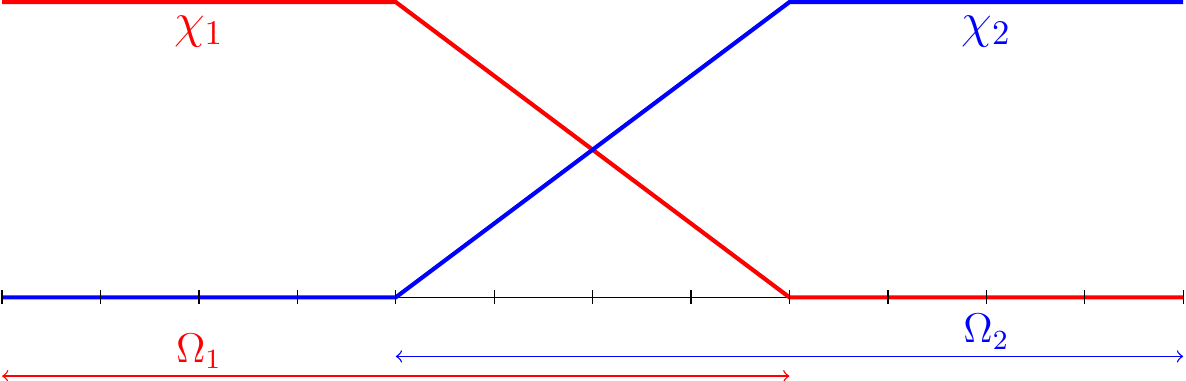}}\\
\subfloat[][{PU1, $\delta=6h$}]
    {\includegraphics[width=.48\textwidth]{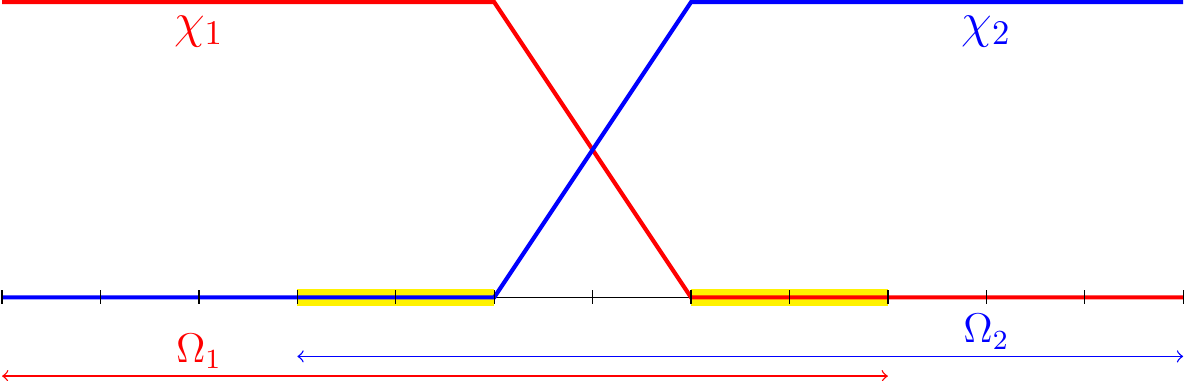}} \quad
\subfloat[][{PU2, $\delta=6h$}]
    {\includegraphics[width=.48\textwidth]{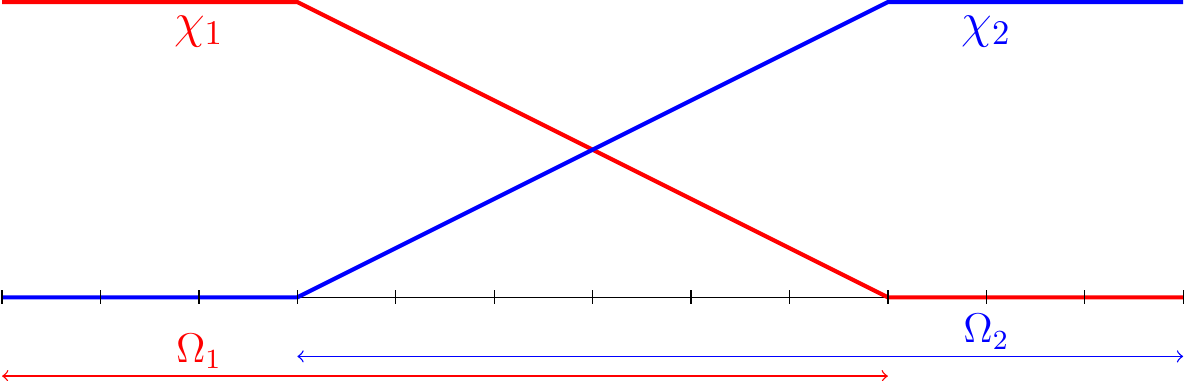}}
\caption{Illustration in a one-dimensional two-subdomain case of the two kinds of partition of unity functions $\chi_j \colon \Omega \to [0,1]$ (PU1 on the left and PU2 on the right), with increasing width of the overlap $\delta$ from top to bottom.}
\label{fig:PU}
\end{figure}
This first choice of Partition of Unity (PU1), where the gradient of $\chi_j$ is zero on the subdomain interfaces $\partial \Omega_j \setminus \partial \Omega$, is illustrated in a one-dimensional two-subdomain case in Figure~\ref{fig:PU}, left, and starting from an overlap $\delta=4h$. Note that PU1 in Figure~\ref{fig:PU} is actually different from the original RAS/ORAS partition of unity, which is defined for any overlap size $\delta$ multiple of $h$, but essentially just at the discrete level, and takes only the values $0$ or $1$; in the original RAS/ORAS articles, the $D_j$ are indeed hidden inside the definition of special extension matrices $\widetilde{R}^T_j$ related to an auxiliary non-overlapping partition of the domain (see e.g.~\cite{EfGa:2003:RAS,StGaTh:2007:OMA} and references therein). However, since the PU1 functions $\chi_j$ in Figure~\ref{fig:PU} are symmetrical to each other, defining the $D_j$ by interpolation of the $\chi_j$ is more practical for a parallel implementation. 

A second kind of Partition of Unity (PU2) is illustrated in Figure~\ref{fig:PU}, right, where the $\chi_j$ functions are different from zero in the interior of the whole overlapping region. This choice is motivated by the fact that using PU1 for SORAS preconditioner can hinder the communication of information between subdomains since in \eqref{eq:SORAS} the matrix $D_j$ is also applied before $B_j^{-1}$, that is before the local problem solve. Indeed, the numerical experiments performed in \cite{BCNT:2021:SORAS}, where PU1 was used, show that the number of iterations for convergence of preconditioned GMRES does not vary significantly when increasing the overlap size (see also Tables~\ref{tab:rotating},\ref{tab:negdiv},\ref{tab:normal} in Section~\ref{sec:num}).

\section{Definition of the model problem} 

As in the second part of \cite{BCNT:2021:SORAS}, we consider the heterogeneous reaction-convection-diffusion problem in conservative form: 
\begin{equation}
\label{eq:RCDbvp}
\begin{cases}
c_0 u + \dive(\mathbf{a}u) -\dive(\nu \nabla u) = f & \text{in } \Omega,\\ 
%\nu \frac{\partial u}{\partial n} -\frac{1}{2} \mathbf{a}\cdot \mathbf{n} \, u + \alpha u = g & \text{on } \Gamma_R,\\
u = 0 & \text{on } \Gamma, %\Gamma_D,
\end{cases}
\end{equation}
where $\Omega \subset \mathbb{R}^d$ is an open bounded polyhedral domain, $\Gamma = \partial \Omega$, 
$\mathbf{n}$ is the outward-pointing unit normal vector to $\Gamma$, $c_0 \in \mathrm{L}^\infty(\Omega)$, $\mathbf{a} \in \mathrm{L}^\infty(\Omega)^d$, $\dive \mathbf{a} \in \mathrm{L}^\infty(\Omega)$, $\nu \in \mathrm{L}^\infty(\Omega)$, $f \in \mathrm{L}^2(\Omega)$ and all quantities are real-valued.  
We denote 
$
\tilde{c} \coloneqq  c_0 + \dive \mathbf{a}/2, 
$
and suppose that there exist $\tilde{c}_- > 0$, $\tilde{c}_+ > 0$ such that 
\begin{equation}
\label{eq:hypctilde}
\tilde{c}_- \le \tilde{c}(\mathbf{x}) \le \tilde{c}_+ \; \text{a.e. in } \Omega,
\end{equation}
%(which is a classical assumption in reaction-convection-diffusion equation literature),
and that there exist $\nu_- > 0$, $\nu_+ > 0$ such that 
$
\nu_- \le \nu(\mathbf{x}) \le \nu_+ \; \text{a.e. in } \Omega.
$ 
The variational formulation of problem~\eqref{eq:RCDbvp} is (see e.g.~\cite[\S4]{BCNT:2021:SORAS}): find $u \in \mathrm{H}_{0}^1(\Omega)$ such that
\begin{equation}
\label{eq:RCDvarf}
a(u,v) = F(v), \quad \text{for all } v \in \mathrm{H}_{0}^1(\Omega),
\end{equation}
%where $a$ is the following non-symmetric bilinear form
\[
a(u,v) \coloneqq \int_\Omega \Bigl ( \tilde{c} u v + \frac{1}{2}  \mathbf{a} \cdot \nabla u \, v - \frac{1}{2} u \, \mathbf{a} \cdot \nabla v  + \nu \nabla u \cdot \nabla v \Bigr ),  \quad F(v) \coloneqq  \int_\Omega fv.
\]
% and 
% \[
% F(v) \coloneqq  \int_\Omega fv.
% \]
On each subdomain we consider the local problem with bilinear form
\[
a_j(u,v) \coloneqq \int_{\Omega_j} \Bigl ( \tilde{c} u v + \frac{1}{2}  \mathbf{a} \cdot \nabla u \, v - \frac{1}{2} u \, \mathbf{a} \cdot \nabla v  + \nu \nabla u \cdot \nabla v \Bigr ) + \int_{\partial \Omega_j \setminus \Gamma} \alpha u v,
\]
where we impose an absorbing transmission condition on the subdomain interface $\partial\Omega_j \setminus \partial\Omega$ given by  
$
\alpha(\mathbf{x}) = \sqrt{(\mathbf{a}\cdot\mathbf{n})^2+4c_0\nu}/2 
$
(see e.g.~\cite{JaNaRo:2001:OO2}).

\section{Numerical experiments}
\label{sec:num}

\begin{table}%[p]
% \caption{Iteration numbers for SORAS preconditioner with the two kinds of partition of unity, in the case of a convection field $\mathbf{a}= 2\pi[-(y-0.1), (x-0.5)]^T$, for different values of the overlap $\delta$, the reaction coefficient $c_0$ and the viscosity $\nu$. 
% The domain is decomposed into $N=5$ overlapping vertical strips and the global problem has 18361 degrees of freedom.} 
\caption{Iteration numbers for SORAS preconditioner ($N=5$).} 
\label{tab:rotating} % aVecCase = 2 in conv-diff-2d-simple_SORAS.edp
\begin{center}
\begin{tabular}{p{4.5cm}*{4}{p{1.1cm}}}
\hline\noalign{\smallskip}
& \multicolumn{4}{c}{\#PU1(PU2)}\\
\cline{2-5}\noalign{\smallskip} $\mathbf{a}= 2\pi[-(y-0.1), (x-0.5)]^T$ & $\delta=2h$ & $\delta=4h$ & $\delta=6h$ & $\delta=8h$ \\
\noalign{\smallskip}\svhline\noalign{\smallskip}
$c_0=1, \; \nu=1$         & 21(21) & 20(17) & 20(15) & 19(14) \\
$c_0=1, \; \nu=0.001$     & 14(14) & 13(11) & 12(11) & 12(10) \\
$c_0=0.001, \; \nu=1$     & 21(21) & 20(18) & 20(15) & 19(14) \\
$c_0=0.001, \; \nu=0.001$ & 15(15) & 14(12) & 13(11) & 13(11) \\
\noalign{\smallskip}\hline\noalign{\smallskip}
\end{tabular}
\end{center}
\end{table}

\begin{table}%[p]
\caption{Repeat of Table~\ref{tab:rotating} but with $\mathbf{a}= [-x, -y]^T$. In this case $\dive \mathbf{a}=-2$ is negative and $\tilde{c}=c_0-1$ does not verify condition~\eqref{eq:hypctilde}.} 
\label{tab:negdiv} % aVecCase = 3 in conv-diff-2d-simple_SORAS.edp
\begin{center}
\begin{tabular}{p{4.5cm}*{4}{p{1.1cm}}}
\hline\noalign{\smallskip}
& \multicolumn{4}{c}{\#PU1(PU2)}\\
\cline{2-5}\noalign{\smallskip} $\mathbf{a}= [-x, -y]^T$ & $\delta=2h$ & $\delta=4h$ & $\delta=6h$ & $\delta=8h$ \\
\noalign{\smallskip}\svhline\noalign{\smallskip}
$c_0=1, \; \nu=1$         & 21(21) & 21(19) & 20(17) & 20(15) \\
$c_0=1, \; \nu=0.001$     & 16(16) & 16(14) & 16(13) & 16(13) \\
$c_0=0.001, \; \nu=1$     & 22(22) & 22(19) & 22(17) & 21(16) \\
$c_0=0.001, \; \nu=0.001$ & 17(17) & 16(15) & 16(14) & 16(13) \\
\noalign{\smallskip}\hline\noalign{\smallskip}
\end{tabular}
\end{center}
\end{table}

\begin{table}%[p]
\caption{Repeat of Table~\ref{tab:rotating} but with $\mathbf{a}=[1, 0]^T$ and with Streamline Upwind Petrov-Galerkin stabilization for the Galerkin approximation.} 
\label{tab:normal} % aVecCase = 1 in conv-diff-2d-simple_SORAS.edp
\begin{center}
\begin{tabular}{p{4.5cm}*{4}{p{1.1cm}}}
\hline\noalign{\smallskip}
& \multicolumn{4}{c}{\#PU1(PU2)}\\
\cline{2-5}\noalign{\smallskip} $\mathbf{a}=[1, 0]^T$ & $\delta=2h$ & $\delta=4h$ & $\delta=6h$ & $\delta=8h$ \\
\noalign{\smallskip}\svhline\noalign{\smallskip}
$c_0=1, \; \nu=1$         & 20(20) & 20(18) & 20(16) & 20(15) \\
$c_0=1, \; \nu=0.001$     & 11(11) & 11(12) & 11(12) & 11(12) \\
$c_0=0.001, \; \nu=1$     & 20(20) & 20(18) & 20(16) & 20(15) \\
$c_0=0.001, \; \nu=0.001$ & 12(12) & 12(12) & 12(13) & 12(12) \\
\noalign{\smallskip}\hline\noalign{\smallskip}
\end{tabular}
\end{center}
\end{table}

% RANDOM INITIAL GUESS
\begin{table}%[p]
% \caption{Weak scaling test, with a regular decomposition into $N$ vertical strips ($\delta=4h$). The global problem has $7381$, $14701$, $29341$, $58621$, $117181$, $234301$ degrees of freedom for $N=2,4,8,16,32,64$ subdomains respectively.} 
\caption{Iteration numbers in a weak scaling test ($\delta=4h$).} 
\label{tab:wscalingReg} % aVecCase = 1  in conv-diff-2d-simple_SORAS.edp
\begin{center}
\begin{tabular}{p{3.5cm}*{6}{p{1.1cm}}}
\hline\noalign{\smallskip}
& \multicolumn{6}{c}{\#PU1(PU2)}\\
\cline{2-7}\noalign{\smallskip}  $\mathbf{a}=[1, 0]^T$ & $N=2$ & $N=4$ & $N=8$ & $N=16$ & $N=32$ & $N=64$ \\
\noalign{\smallskip}\svhline\noalign{\smallskip}
$c_0=1, \; \nu=1$         & 18(15) & 23(20) & 28(24) & 35(28) & 36(29) & 36(29) \\
$c_0=1, \; \nu=0.001$     & 8(8) & 10(12) & 16(16) & 23(24) & 37(37) & 63(61) \\
$c_0=0.001, \; \nu=1$     & 18(15) & 23(20) & 29(25) & 35(29) & 36(29) & 36(29) \\
$c_0=0.001, \; \nu=0.001$ & 8(8) & 10(12) & 16(17) & 24(25) & 40(40) & 71(71) \\
\noalign{\smallskip}\hline\noalign{\smallskip}
\end{tabular}
\end{center}
\end{table}

% \begin{table}%[p]
% \caption{Weak scaling test as in Table~\ref{tab:wscalingReg}, but with METIS decomposition into $N$ arbitrary-shaped subdomains.} 
% \label{tab:wscalingMetis} % aVecCase = 1  in conv-diff-2d-simple_SORAS.edp
% \begin{center}
% \begin{tabular}{p{3.5cm}*{6}{p{1.1cm}}}
% \hline\noalign{\smallskip}
% & \multicolumn{6}{c}{\#PU1(PU2)}\\
% \cline{2-7}\noalign{\smallskip}  $\mathbf{a}=[1, 0]^T$ & $N=2$ & $N=4$ & $N=8$ & $N=16$ & $N=32$ & $N=64$ \\
% \noalign{\smallskip}\svhline\noalign{\smallskip}
% $c_0=1, \; \nu=1$         & 21(18) & 30(26) & 40(33) & 48(39) & 53(41) & 55(43) \\
% $c_0=1, \; \nu=0.001$     & 10(9) & 12(12) & 17(17) & 25(24) & 38(38) & 63(62) \\
% $c_0=0.001, \; \nu=1$     & 21(18) & 30(26) & 40(33) & 48(40) & 54(42) & 57(44) \\
% $c_0=0.001, \; \nu=0.001$ & 10(10) & 12(12) & 18(18) & 26(26) & 42(41) & 73(72) \\
% \noalign{\smallskip}\hline\noalign{\smallskip}
% \end{tabular}
% \end{center}
% \end{table}

We simulate problem~\eqref{eq:RCDvarf} with $\Omega$ a rectangle $[0,N\cdot0.2]\times[0,0.2]$, where $N$ is the number of subdomains. %, and $\Gamma_D = \Gamma$, $\Gamma_R=\emptyset$. 
In Tables~\ref{tab:rotating},\ref{tab:negdiv},\ref{tab:normal} we take $N=5$ and $f=100\exp{\{-10((x-0.5)^2+(y-0.1)^2)\}}$. In Table~\ref{tab:wscalingReg}, we test weak scaling by  
varying $N$, with $f=100\exp{\{-10((x-0.1)^2+(y-0.1)^2)\}}$. 
The problem is discretized by piece-wise linear Lagrange finite elements on a uniform triangular mesh with $60$ {nodes} on the vertical side of the rectangle and $N \cdot 60$ {nodes} on the horizontal one, resulting in $18361$ degrees of freedom for $N=5$, and $7381$, $14701$, $29341$, $58621$, $117181$, $234301$ degrees of freedom for $N=2,4,8,16,32,64$ respectively.  
The domain is partitioned into $N$ vertical strips, then each subdomain is augmented with mesh elements layers of size $\delta/2$ to obtain the overlapping decomposition: the total width of the overlap between two subdomains is then $\delta$. {In particular, for $\delta=2h, 4h, 6h, 8h$ the ratio between the subdomain width ($60h$) and $\delta$ is equal to $30, 15, 10, 7.5$.} 
We use GMRES with right preconditioning, with a zero initial guess in Tables~\ref{tab:rotating},\ref{tab:negdiv},\ref{tab:normal} and a random initial guess in Table~\ref{tab:wscalingReg}. The stopping criterion is based on the relative residual, with a tolerance of $10^{-6}$. 
To apply the preconditioner, the local problems in each subdomain are solved with the direct solver \href{http://mumps.enseeiht.fr/}{MUMPS}.  
All the computations are done in the \texttt{ffddm} framework \cite{FFD:Tournier:2020} of \href{https://freefem.org/}{FreeFEM}. 

We compare the number of iterations for convergence (denoted by \# in the tables) using the two kinds of partition of unity: the results for PU1 were also included in \cite{BCNT:2021:SORAS} and the results for PU2 are reported inside brackets in Tables~\ref{tab:rotating}--\ref{tab:wscalingReg}. In \texttt{ffddm} framework, the first partition of unity is selected by the flag \texttt{-raspart}, while the second partition of unity is the one used by default.   

As in \cite{BCNT:2021:SORAS}, we examine several configurations for the coefficients in \eqref{eq:RCDbvp}. First, in Table~\ref{tab:rotating} we consider a rotating convection field $\mathbf{a}= 2\pi[-(y-0.1), (x-0.5)]^T$ and small/large values for the reaction coefficient $c_0$ and the viscosity $\nu$. We can see that a larger overlap helps the convergence of the preconditioner, especially with PU2, while with PU1 the number of iterations does not vary significantly. Moreover, with both kinds of partition of unity, the number of iterations appears not very sensitive to the reaction coefficient $c_0$, while it increases when the viscosity $\nu$ is larger. 

Then, in Table~\ref{tab:negdiv} we take $\mathbf{a}= [-x, -y]^T$, which has negative divergence $\dive \mathbf{a}=-2$, to test the robustness of the method when condition~\eqref{eq:hypctilde} on the positiveness of $\tilde{c}$ is violated: in this case, $\tilde{c}=c_0-1$, so $\tilde{c}=0$,  $\tilde{c}=-0.999$ for $c_0=1$, $c_0=0.001$ respectively. We can still observe a convergence behavior similar to the one of Table~\ref{tab:rotating}. 

Finally, in Table~\ref{tab:normal} we consider a horizontal convection field $\mathbf{a}=[1, 0]^T$, which is normal to the interfaces between subdomains. Since in this case non-physical numerical instabilities appear in the solution, we stabilize the discrete variational formulation using the Streamline Upwind Petrov-Galerkin (SUPG) method (see for instance \cite[\S11.8.6]{quarteroni:2009:book}).
%, which adds to the Galerkin approximation the following term (see for instance \cite[\S11.8.6]{quarteroni:2009:book}): 
% \[
% \mathcal{L}_h(u_h,f;v_h) = \theta \sum_{\tau\in \mathcal{T}^h} \int_\tau (\mathcal{L}u_h-f) \, \frac{h_\tau}{\rvert \mathbf{a}\rvert} \, \mathcal{L}_{SS}v_h, 
% \]
% where $\theta$ is a stabilization parameter (here we choose $\theta=0.15$), $h_\tau$ is the diameter of the mesh element $\tau$, and 
% \[
% \mathcal{L}u_h = c_0 u_h + \dive(\mathbf{a}u_h) -\dive(\nu \nabla u_h), \quad
% \mathcal{L}_{SS}v_h = \frac{1}{2}\dive(\mathbf{a}v_h)+\frac{1}{2}\mathbf{a}\cdot\nabla v_h.
% \]
In this configuration for the convection field, for low viscosity $\nu=0.001$ the dependence of the iteration number on the overlap size $\delta$ appears to be not significant, even with PU2. 

Again in this third configuration with $\mathbf{a}=[1, 0]^T$ and SUPG stabilization, we perform a weak scaling test by taking $\Omega = [0,N\cdot0.2]\times[0,0.2]$ for increasing number of subdomains $N$, and $\delta=4h$.  
%First we consider a regular decomposition into vertical strips (Table~\ref{tab:wscalingReg}) and then an arbitrary decomposition obtained by METIS (Table~\ref{tab:wscalingMetis}). In both tables, we fix the overlap $\delta=4h$. 
% Comparing Table~\ref{tab:wscalingReg} with Table~\ref{tab:wscalingMetis}, we can see that the number of iterations is higher when taking arbitrary-shaped subdomains. Moreover, 
We can see that especially in the cases with low viscosity $\nu=0.001$, convergence deteriorates with $N$, as expected since we are testing a one-level preconditioner. 
%In most tests, using PU2 improves the iteration counts obtained with PU1. 

In summary, our numerical investigation shows that, for the considered SORAS preconditioner, PU2 generally improves the iteration counts obtained with PU1.  
Moreover, the first kind of partition of unity (PU1), which would be the natural choice for ORAS solver instead, yields for SORAS preconditioner iterations counts that do not vary significantly when increasing the overlap width, {whereas using the second kind of partition of unity (PU2) a larger overlap gives faster convergence.}

\begin{table}%[p]
\caption{Minimum and maximum eigenvalues of the preconditioned operator.} 
\label{tab:spectrasym} % aVecCase = 3 in conv-diff-2d-simple_SORAS.edp
\begin{center}
\begin{tabular}{p{1.8cm}*{4}{p{1.8cm}}}
\hline\noalign{\smallskip}
& \multicolumn{4}{c}{PU1(PU2)}\\
\cline{2-5}\noalign{\smallskip} $\mathbf{a}= [0, 0]^T$ & $\delta=2h$ & $\delta=4h$ & $\delta=6h$ & $\delta=8h$ \\
\noalign{\smallskip}\svhline\noalign{\smallskip}
$\lambda_{\text{min}}$ & 0.50 (0.50) & 0.50 (0.50) & 0.50 (0.50) & 0.50 (0.50)\\
$\lambda_{\text{max}}$ & 11.25 (11.25) & 10.61 (5.98) & 10.07 (4.01) & 9.60 (3.02)\\
\noalign{\smallskip}\hline\noalign{\smallskip}
\end{tabular}
\end{center}
\end{table}

To conclude, we wish to provide a deeper explanation of the observed effects.  
First, we examine the symmetric positive definite case, with $\mathbf{a}= [0, 0]^T$, $c_0 = 1$, $\nu=1$ (so $\tilde{c}=c_0>0$), and report in Table~\ref{tab:spectrasym} the largest and smallest eigenvalues of the preconditioned operator. We take $N=2$ and $40$ nodes on the vertical side of the rectangle, $2 \cdot 40$ nodes on the horizontal one. Note that SORAS preconditioner for generic symmetric positive definite problems was analyzed in \cite{HaJoNa:soras:2015}, but no explicit discussion about the influence of the partition of unity was included there. 
On the one hand, the largest eigenvalue of the preconditioned operator is controlled by the modes of the local generalized eigenvalue problems defined in \cite[Definition 3.1]{HaJoNa:soras:2015}, where the partition of unity matrices appear in the local operator on the left-hand side: Table~\ref{tab:spectrasym} shows that indeed $\lambda_{\text{max}}$ is smaller for PU2, which is less steep than PU1, especially when increasing the overlap width $\delta$ (see Fig.~1). Moreover, with PU1, the dependence of $\lambda_{\text{max}}$ on $\delta$ is much less significant than with PU2.   
On the other hand, the smallest eigenvalue of the preconditioned operator is controlled by the modes of the local generalized eigenvalue problems defined in \cite[Definition 3.2]{HaJoNa:soras:2015}, where the partition of unity is not involved: in Table~\ref{tab:spectrasym} we can see that indeed $\lambda_{\text{min}}$ is independent of the partition of unity.

\begin{figure}
    \centering
    \subfloat[][PU1]
        {\includegraphics[width=.5\textwidth]{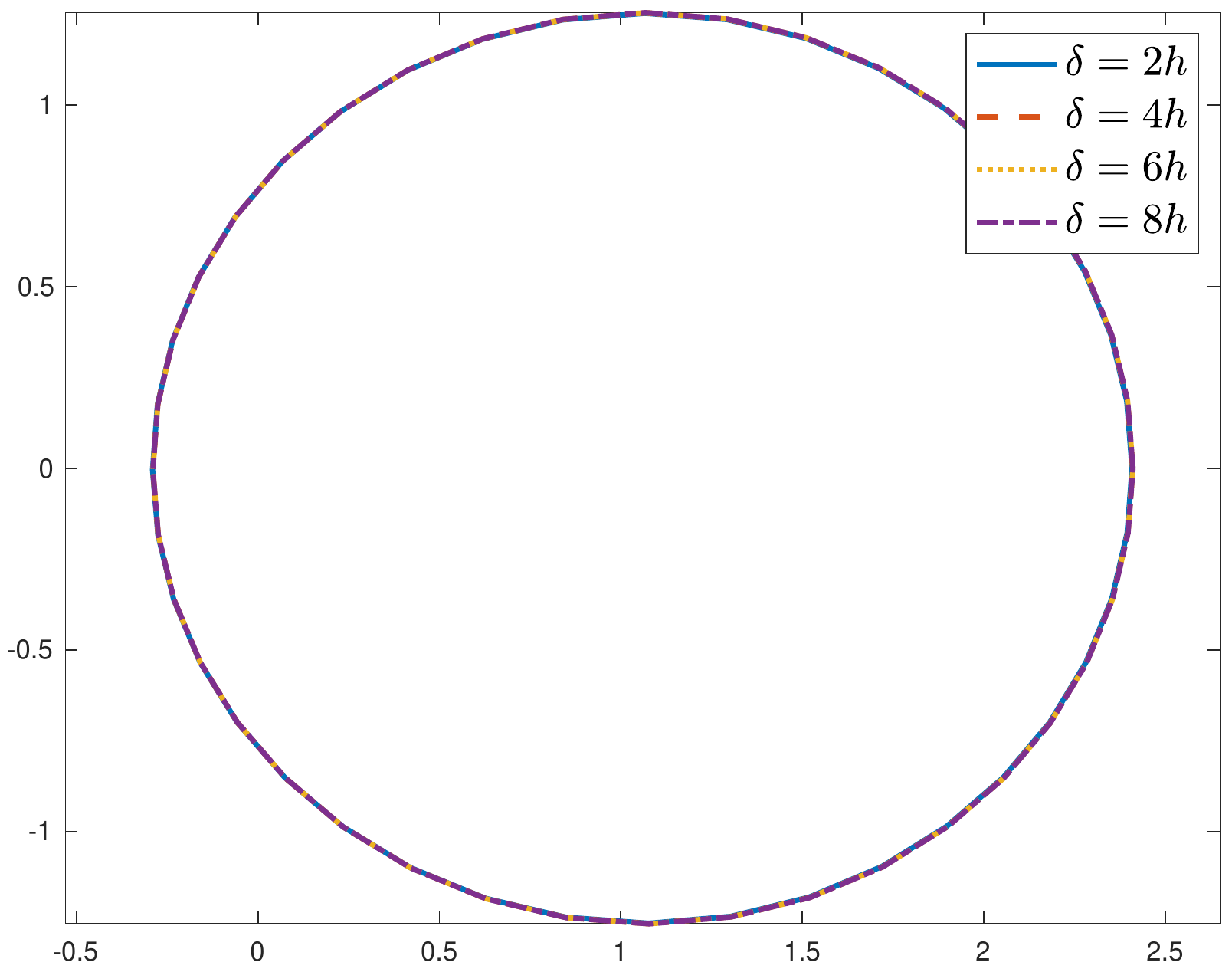}} 
    \subfloat[][PU2]
        {\includegraphics[width=.5\textwidth]{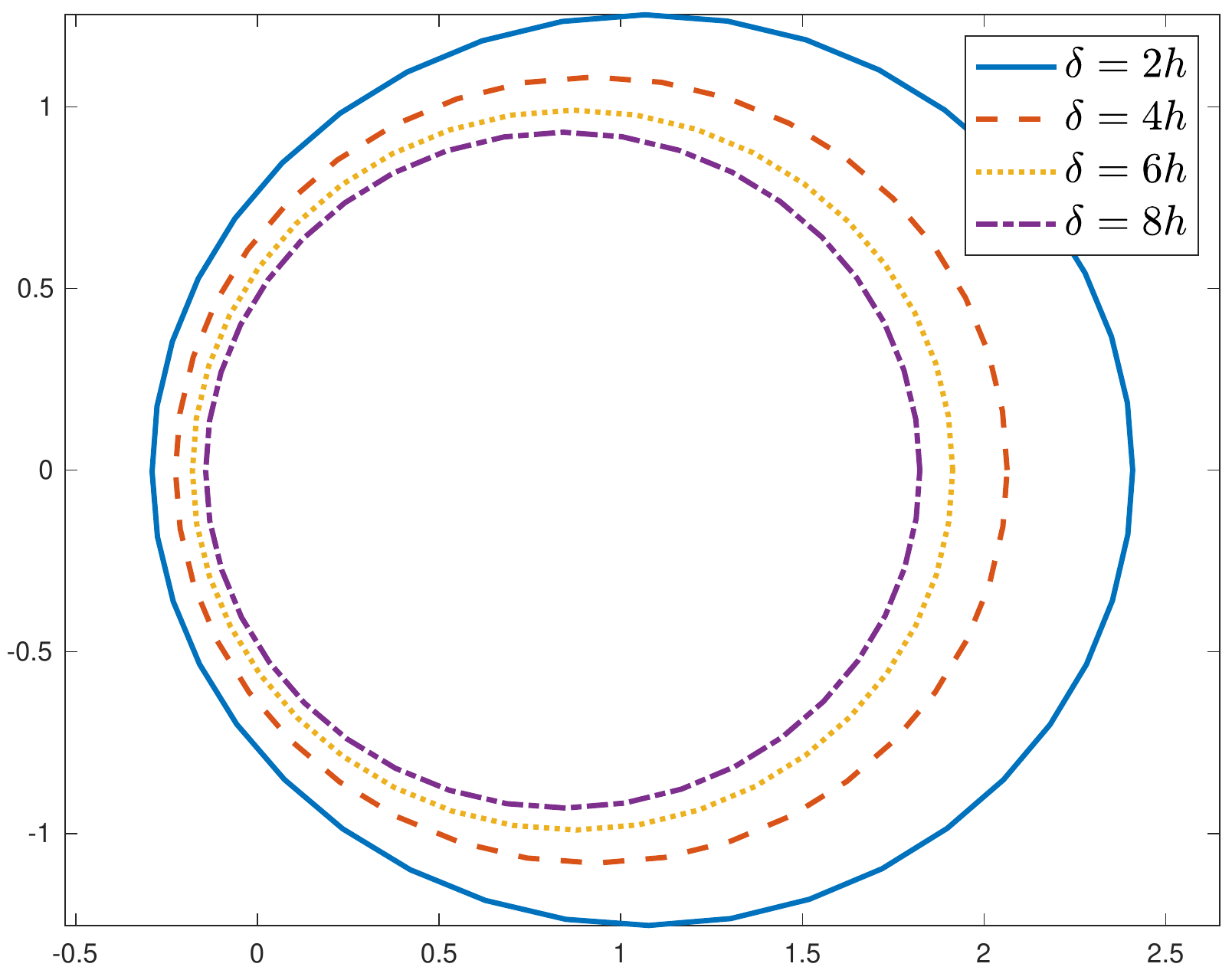}}
    \caption{Numerical range of the preconditioned operator ($\mathbf{a}= [-x, -y]^T$).}
    \label{fig:fov}
\end{figure}

For the non-symmetric case, with $\mathbf{a}= [-x, -y]^T$, $c_0=0.001$, $\nu=0.001$ (so $\tilde{c}=c_0-1<0$), we plot in Fig.~\ref{fig:fov} the contour of the numerical range of the preconditioned operator for overlap widths that range from $\delta=2h$ to $\delta=8h$, for the two types of partition of unity. We can remark that for PU1 (Fig.~\ref{fig:fov}, left) the numerical ranges practically coincide for the different overlap widths, whereas for PU2 (Fig.~\ref{fig:fov}, right) the numerical range gets smaller for larger overlap width. This explains the more favorable convergence properties of preconditioned GMRES with PU2 when increasing the overlap width, and the much less significant influence of the overlap in the case of PU1.     

\bibliographystyle{spmpsci}
\bibliography{authorsample}

\end{document}